\documentclass[12pt]{article} % default is 10 pt

\usepackage{latexsym,amsmath,amssymb,a4wide,amscd}
\usepackage[latin1]{inputenc}
\usepackage[T1]{fontenc}
\newcommand{\changefont}[3]{
\fontfamily{#1} \fontseries{#2} \fontshape{#3} \selectfont}

\usepackage{array}

\usepackage{amsthm}

\usepackage{url}

\usepackage{mathabx}

\usepackage{longtable}

\numberwithin{equation}{section}
\numberwithin{figure}{section}

\swapnumbers
\theoremstyle{definition}

\usepackage[pdftex]{graphicx} % needed for including graphics e.g. EPS, PS
\usepackage[nooneline]{caption}

\begin{document}
% don't want date printed
\date{}

\title{\Large {\bf Supplemental material to the article ``Partitions of the triangles of the cross polytope into surfaces''}}

\author{Jonathan Spreer}

\maketitle

\subsection*{\centering Abstract}

{\em
	We present a constructive proof, that there exists a decomposition of the $2$-skeleton of the $k$-dimensional cross polytope $\beta^k$ into closed surfaces of genus $\leq 1$, each with a transitive automorphism group given by the vertex transitive $\mathbb{Z}_{2k}$-action on $\beta^k$. Furthermore we show, that for each $k \equiv 1,5(6)$ the $2$-skeleton of the $(k-1)$-simplex is a union of highly symmetric tori and M\"obius strips.
}\\
\\
\textbf{MSC 2010: } {\bf 52B12}; % centrally symmetric polytopes
52B70; % polytopal manifolds
57Q15;  % triangulating manifolds
57M20; % Two-dimensional complexes  
05C10; % Planar graphs; geometric and topological aspects of graph theory
\\
\textbf{Keywords: } cross polytope, simplicial complexes, triangulated surfaces, difference cycles.

\vspace{1.0cm}
The following lists contain detailed information about the decompositions of $\beta^k$ and $\Delta^{k-1}$ described in the Theorems $2.1$ and $3.1$ of the article ``Partitions of the triangles of the cross polytope into surfaces''.

\scriptsize
\begin{center}
\begin{longtable}{|l|l|l@{}l@{}l|}
		\caption{The decomposition of the $2$-skeleton of $\beta^{k}$ ($k \leq 90 $) by topological types.} \label{tab:listBeta} \\
		
		\hline \multicolumn{1}{|l|}{$k$} & \multicolumn{1}{l|}{\textbf{topological type}} & \multicolumn{3}{l|}{\textbf{difference cycles}} \\ \hline 
		\endfirsthead

		\multicolumn{5}{l}%
		{ \tablename\ \thetable{} -- continued from previous page} \\
		\hline \multicolumn{1}{|l|}{$k$} & \multicolumn{1}{l|}{\textbf{topological type}} & \multicolumn{3}{l|}{\textbf{difference cycles}} \\ \hline 
		\endhead

		\hline \multicolumn{5}{r}{{continued on next page --}} \\
		\endfoot

		\hline \hline
		\endlastfoot

$3$&$\mathbb{S}^2$&$\{(2\!:\!2\!:\!2), (1\!:\!1\!:\!4)\}$&& \\
\hline 
\hline 
$4$&$\mathbb{T}^2$&$\{ ( 1\!:\!2\!:\!5 ), ( 1\!:\!5\!:\!2 ) \},$&$\{ ( 1\!:\!1\!:\!6 ), ( 3\!:\!3\!:\!2 ) \}$&\\ 
\hline 
\hline 
$5$&$\mathbb{T}^2$&$\{ ( 1\!:\!2\!:\!7 ), ( 1\!:\!7\!:\!2 ) \},$&$\{ ( 1\!:\!3\!:\!6 ), ( 1\!:\!6\!:\!3 ) \}$&\\ 
\hline 
&$\mathbb{K}^2$&$\{ ( 1\!:\!1\!:\!8 ), ( 4\!:\!4\!:\!2 ) \},$&$\{ ( 2\!:\!2\!:\!6 ), ( 3\!:\!3\!:\!4 ) \}$&\\ 
\hline 
\hline 
$6$&$\{ 1 , 2 \} \times \mathbb{S}^2$&$\{(4\!:\!4\!:\!4), (2\!:\!2\!:\!8)\}$&& \\
\hline 
&$\mathbb{T}^2$&$\{ ( 1\!:\!2\!:\!9 ), ( 1\!:\!9\!:\!2 ) \},$&$\{ ( 1\!:\!3\!:\!8 ), ( 1\!:\!8\!:\!3 ) \},$&$\{ ( 1\!:\!4\!:\!7 ), ( 1\!:\!7\!:\!4 ) \},$\\ 
&&$\{ ( 2\!:\!3\!:\!7 ), ( 2\!:\!7\!:\!3 ) \},$&$\{ ( 3\!:\!4\!:\!5 ), ( 3\!:\!5\!:\!4 ) \},$&$\{ ( 1\!:\!1\!:\!10 ), ( 5\!:\!5\!:\!2 ) \}$\\ 
\hline 
\hline 
$7$&$\mathbb{T}^2$&$\{ ( 1\!:\!2\!:\!11 ), ( 1\!:\!11\!:\!2 ) \},$&$\{ ( 1\!:\!3\!:\!10 ), ( 1\!:\!10\!:\!3 ) \},$&$\{ ( 1\!:\!4\!:\!9 ), ( 1\!:\!9\!:\!4 ) \},$\\
\hline 
\hline 
		\end{longtable}
\end{center}

\scriptsize

\begin{center}
	\begin{longtable}{|l|l|l@{}l@{}l|}
		\caption{The decomposition of the $2$-skeleton of $\Delta^{k-1}$ ($k \leq 100 $) by topological types.} \label{tab:listDelta} \\
		
		\hline \multicolumn{1}{|l|}{$k$} & \multicolumn{1}{l|}{\textbf{topological type}} & \multicolumn{3}{l|}{\textbf{difference cycles}} \\ \hline 
		\endfirsthead

		\multicolumn{5}{l}%
		{ \tablename\ \thetable{} -- continued from previous page} \\
		\hline \multicolumn{1}{|l|}{$k$} & \multicolumn{1}{l|}{\textbf{topological type}} & \multicolumn{3}{l|}{\textbf{difference cycles}} \\ \hline 
		\endhead

		\hline \multicolumn{5}{r}{{continued on next page --}} \\
		\endfoot

		\hline \hline
		\endlastfoot

$5$&$\mathbb{M}^2$&$\{ ( 1\!:\!1\!:\!3 ) \},$&$\{ ( 2\!:\!2\!:\!1 ) \}$&\\ 
\hline 
\hline 
$7$&$\mathbb{M}^2$&$\{ ( 1\!:\!1\!:\!5 ) \},$&$\{ ( 2\!:\!2\!:\!3 ) \},$&$\{ ( 3\!:\!3\!:\!1 ) \}$\\ 
\hline 
&$\mathbb{T}^2$&$\{ ( 1\!:\!2\!:\!4 ), ( 1\!:\!4\!:\!2 ) \}$&&\\ 
\hline 
\hline 
$11$&$\mathbb{M}^2$&$\{ ( 1\!:\!1\!:\!9 ) \},$&$\{ ( 2\!:\!2\!:\!7 ) \},$&$\{ ( 3\!:\!3\!:\!5 ) \},$\\ 
&&$\{ ( 4\!:\!4\!:\!3 ) \},$&$\{ ( 5\!:\!5\!:\!1 ) \}$&\\ 
\hline 
&$\mathbb{T}^2$&$\{ ( 1\!:\!2\!:\!8 ), ( 1\!:\!8\!:\!2 ) \},$&$\{ ( 1\!:\!3\!:\!7 ), ( 1\!:\!7\!:\!3 ) \},$&$\{ ( 1\!:\!4\!:\!6 ), ( 1\!:\!6\!:\!4 ) \},$\\ 
&&$\{ ( 2\!:\!3\!:\!6 ), ( 2\!:\!6\!:\!3 ) \},$&$\{ ( 2\!:\!4\!:\!5 ), ( 2\!:\!5\!:\!4 ) \}$&\\ 
\hline 
\hline 
$13$&$\mathbb{M}^2$&$\{ ( 1\!:\!1\!:\!11 ) \},$&$\{ ( 2\!:\!2\!:\!9 ) \},$&$\{ ( 3\!:\!3\!:\!7 ) \},$\\
\hline 
\hline 
		\end{longtable}
\end{center}

\noindent
Jonathan Spreer \\
Institut f\"ur Geometrie und Topologie \\
Universit\"at Stuttgart \\
70550 Stuttgart \\
Germany
\end{document}